\documentclass[12pt]{amsart}

\usepackage{amsmath,amssymb}
\begin{document}
\newtheorem{thm}{Theorem}[section]
\newtheorem{lem}[thm]{Lemma}
\newtheorem{dfn}[thm]{Definition}
\newtheorem{cor}[thm]{Corollary}
\newtheorem{conj}[thm]{Conjecture}
\newtheorem{clm}[thm]{Claim}
\theoremstyle{remark}
\newtheorem{exm}[thm]{Example}
\newtheorem{rem}[thm]{Remark}
\def\N{{\mathbb N}}
\def\G{{\mathbb G}}
\def\Q{{\mathbb Q}}
\def\R{{\mathbb R}}
\def\C{{\mathbb C}}
\def\P{{\mathbb P}}
\def\F{{\mathbb F}}
\def\Z{{\mathbb Z}}
\def\v{{\mathbf v}}
\def\x{{\mathbf x}}
\def\O{{\mathcal O}}
\def\M{{\mathcal M}}
\def\PP{{\mathcal P}}
\def\tr{\mbox{Tr}}
\def\ord{{\mbox{ord}}}
\def\qed{{\tiny $\clubsuit$ \normalsize}}
\def\dist{{\mbox{dist}}}

\title[Vojta implies Batyrev-Manin]{Vojta's Conjecture implies the
  Batyrev-Manin Conjecture for $K3$ surfaces}

\author{David McKinnon}
\thanks{Supported in part by a Discovery Grant from the Natural Sciences and
Engineering Research Council of Canada}

\begin{abstract}
Vojta's Conjectures are well known to imply a wide range of results,
known and unknown, in arithmetic geometry.  In this paper, we add to
the list by proving that they imply that rational points tend to repel
each other on algebraic varieties with nonnegative Kodaira dimension.
We use this to deduce, from Vojta's Conjectures, conjectures of
Batyrev-Manin and Manin on the distribution of rational points on
algebraic varieties.  In particular, we show that Vojta's Main
Conjecture implies the Batyrev-Manin Conjecture for $K3$ surfaces.
\end{abstract}

\maketitle

\section{Introduction}

In 1987, in \cite{Vo}, Vojta made a series of wide-ranging and deep
conjectures about the distribution of rational points on algebraic
varieties.  They imply the Masser-Oesterl\'e $abc$ conjecture, the
Bombieri-Lang Conjecture that the set of rational points on a variety
of general type is not Zariski dense, and a host of other well known
conjectures in number theory.  (For more details on the various
implications of Vojta's conjectures, see \cite{Vo} or \cite{HS},
section F.5.3.)

Vojta's Main Conjecture is known in only a few special cases, although
some of these special cases are extremely significant.  It is known
for curves, being the union of Roth's Theorem for rational curves
(\cite{Ro}), Siegel's Theorem for elliptic curves (\cite{Sg}), and
Faltings' Theorem for curves of general type (\cite{Fa1}).  It is also
known for abelian varieties, by another theorem of Faltings
(\cite{Fa2}).  The case $X=\P^n$ and $D$ a union of hyperplanes is
also known, due originally to Schmidt for archimedean places $S$
(\cite{Sm}), and Schlickewei for a general set of places (\cite{Sl}),
and is the famous Schmidt Subspace Theorem.  For a precise statement of
Vojta's Main Conjecture, see section~\ref{statement}.

For varieties with a canonical class of which some multiple is
effective (that is, for varieties with nonnegative Kodaira dimension),
there is a lower bound for a height with respect to the canonical
class.  In this case, Vojta's conjectures imply that rational points
repel each other, in the sense that two rational points that are close
to one another (with respect to, say, some archimedean metric) must
have large height relative to this distance.  This idea is made
precise by the Repulsion Principle, Theorem~\ref{technical}, in
section~\ref{repulsion}.

From the Repulsion Principle, one can easily deduce that on an
algebraic variety of nonnegative Kodaira dimension, rational points
must be quite sparse.  For $K3$ surfaces, this is precisely the famous
conjecture of Batyrev and Manin (\cite{BM}), and in section~\ref{bm}
we deduce this conjecture from Vojta's Main Conjecture, as a corollary
of a more general result (Theorem~\ref{main}).  Finally, in
section~\ref{rcc}, we deduce Manin's Rational Curve Conjecture
(\cite{Ma}) from the Repulsion Principle.

It is a pleasure to thank Ekaterina Amerik for her helpful remarks
that helped to simplify and clarify the arguments in this paper.

\section{Vojta's Main Conjecture}\label{statement}

Vojta's Main Conjecture requires a lot of terminology to state, so we
list some of it here for convenience.

\vspace{.1in}

\begin{tabular}{cl}
$k$ & A number field \\
$S$ & A finite set of places of $k$ \\
$X$ & A smooth algebraic variety defined over $k$ \\
$K$ & The canonical divisor class of $X$ \\
$D$ & A normal crossings divisor on $X$ 
(see Definition~\ref{normalcrossings}) \\
$L$ & A big divisor on $X$ (see Definition~\ref{big}) \\
$h_K$, $h_L$ & Logarithmic height functions with respect to $K$ and $L$ \\
$h_{D,v}$ & A local height function for $D$ with respect to a place 
$v$ \\
& of $k$ (see Definition~\ref{localheight}) \\
$m_S(D,\cdot)$ & A proximity function for $D$ with respect to $S$, given
by \\
& $m_S(D,P)=\sum_{v\in S} h_{D,v}(P)$
\end{tabular}

\vspace{.1in}

The conjecture we are primarily interested in here is the following
(see Conjecture~3.4.3 from \cite{Vo}):

\begin{conj}[Vojta's Main Conjecture]\label{voriginal}
Choose any $\epsilon>0$.  Then there exists a nonempty Zariski open
set $U=U(\epsilon)\subset X$ such that for every $k$-rational point
$P\in U(k)$, we have the following inequality:
\begin{equation}
m_S(D,P) + h_K(P) \leq \epsilon h_L(P)
\end{equation}
\end{conj}

For a general discussion of this conjecture, and an explanation of all
the terminology, we refer the reader to \cite{Vo}.  However, for the
reader's convenience, we recall the definitions of some of the less
common terms mentioned in Conjecture~\ref{voriginal}:

\begin{dfn}\label{big}
A divisor $D$ on a smooth algebraic variety $X$ is {\bf big} if and
only if it can be written as the sum $D=A+E$ of an ample divisor $A$
and an effective divisor $E$.  Note that in particular every ample
divisor is big.
\end{dfn}

\begin{dfn}\label{normalcrossings}
Let $X$ be a smooth algebraic variety defined over the complex
numbers.  A divisor $D$ on $X$ has {\bf normal crossings} if and only
if it is effective, has no multiple components, and for every point
$P$ in the support of $D$, there are analytic functions
$z_1,\ldots,z_n$ on an analytic neighbourhood $U$ of $P$ such that
$D\cap U$ is the locus ${z_1\ldots z_n=0}\cap U$, where $n\leq\dim X$.

For a smooth algebraic variety $X$ defined over a number field $k$, we
say that a divisor $D$ has {\bf normal crossings} if and only if for
every embedding $k\hookrightarrow\C$, the corresponding divisor $D_\C$
on $X_\C$ has normal crossings.  (It is easy to see that if $D_\C$ has
normal crossings for any one embedding of $k$ in $\C$, then $D_\C$
will have normal crossings for all such embeddings.)
\end{dfn}

\begin{dfn}\label{localheight}
Let $X$ be a smooth algebraic variety defined over a number field $k$,
let $D$ be a divisor on $X$, and let $v$ be a place of $k$.  Let $k_v$
be the completion of $k$ with respect to $v$, and let $\overline{k_v}$
be an algebraic closure of $k_v$.  Let $\|\cdot\|_v$ be the absolute
value corresponding to $v$, extended to an absolute value on
$\overline{k_v}$.  A {\bf local height function} for $D$ at $v$ is a
function $h_{D,v}\colon X(\overline{k_v})-\mbox{Supp}(D)\to\R$ such
that for all $P$ in the support $\mbox{Supp}(D)$ of $D$, there is a
Zariski open neighbourhood $U$ of $P$ such that $D\cap U$ is the locus
${f=0}\cap U$ for some rational function $f$, and:
\[h_{D,v}(P) = -\frac{1}{[k:\Q]}\log\|f(P)\|_v+\alpha(P)\]
where $\alpha$ is a continuous function on $U(\overline{k_v})$.

If $D$ is an effective cycle on $X$ that can be written as the
scheme-theoretic intersection of a finite number of effective divisors
$D=\bigcap_i D_i$, then we define $h_{D,v}=\min_i h_{D_i,v}$.  Note that
$h_{D,v}$ can be defined on $X(\overline{k_v})-\mbox{Supp}(D)$ by
setting $h_{D_i,v}(P)=\infty$ if $P\in\mbox{Supp}(D_i)$.
\end{dfn}

It is not hard to see that if one chooses a different set of divisors
$\{D_i\}$, then the resulting local height function differs from the
original by a bounded function.

Intuitively, one may think of the local height functions as satisfying
$h_{D,v}(P)=-\log\dist_v(P,D)$, where $\dist_v(P,D)$ denotes the
distance from $P$ to the support of $D$.  Note that in \cite{Vo},
local height functions are called Weil functions.

In what follows, we will want to apply Conjecture~\ref{voriginal} to
the slightly more general case in which $D$ is a cycle of arbitrary
codimension.  To do this, we return to and slightly generalise the
notation from Vojta's Main Conjecture (Conjecture~\ref{voriginal}):

\vspace{.1in}

\begin{tabular}{cl}
$D'$ & A cycle on $X$ that is contained in a normal crossings divisor
$D$ \\
$h_{D',v}$ & A local height function for $D'$ with respect to $v$ \\
$m_S(D',\cdot)$ & A proximity function for $D'$ with respect to $s$, 
given by \\
& $m_S(D',P) = \sum_{v\in S} h_{D',v}(P)$
\end{tabular}

\begin{conj}[Vojta's Main Conjecture for general cycles]\label{vojta}
Choose any $\epsilon>0$.  Then there exists a nonempty Zariski open
set $U=U(\epsilon)\subset X$ such that for every $k$-rational point
$P\in U(k)$, we have the following inequality:
\begin{equation}
m_S(D',P) + h_K(P) \leq \epsilon h_L(P)
\end{equation}
\end{conj}

This follows immediately from Vojta's Main Conjecture and the observation
that $m_S(D',P)\leq m_S(D,P)$.

\section{Repulsion of rational points for varieties of nonnegative 
Kodaira dimension}\label{repulsion}

This section contains a technical result which contains the heart of
the proofs of the main results of the paper.  In particular,
Theorem~\ref{technical} implies that the rational points of low height
on a variety of nonnegative Kodaira dimension (such as a $K3$ surface)
should be distributed very sparsely away from subvarieties of negative
Kodaira dimension (such as rational curves).  Specifically, we have:

\begin{thm}[Repulsion Principle]\label{technical}
Let $X$ be any smooth, projective variety of nonnegative Kodaira
dimension defined over a number field $k$.  For every smooth variety
$V$ birational to a subvariety of $X\times X$, make the following two
assumptions:
\begin{enumerate}
\item If $V$ has nonnegative Kodaira dimension, then it satisfies
  Vojta's Main Conjecture for general cycles (Conjecture~\ref{vojta}).
\item If $V$ has negative Kodaira dimension, then it is uniruled.  
\end{enumerate}
Let $v$ be a place of $k$, and choose an ample multiplicative height
$H$ on $Y$.

Then for any $\epsilon>0$, there is a nonempty Zariski open subset
$U(\epsilon)$ of $X$ and a positive real constant $C$ such that
\[\dist_v(P,Q) > CH(P,Q)^{-\epsilon}\]
for all $P,Q\in U(k)$.
\end{thm}

\noindent
{\it Proof:} \/ Let $v$ be a place of $k$.  Let $Y$ be the variety
$X\times X$, and let $D$ be the diagonal on $Y$.  Let $h$ be an ample
logarithmic height on $X$, and define an ample logarithmic height on
$Y$ by $h(P,Q)=h(P)+h(Q)$.  By Vojta's Main Conjecture for general
cycles (Conjecture~\ref{vojta}), there is a proper Zariski closed
subset $Z\subset Y$ such that for every $(P,Q)\in Y(k)-Z(k)$, we have:
\begin{align}
m_{D,v}(P,Q) &\leq \epsilon h(P,Q) + O(1) \label{conjeqn}
\end{align}
In what follows, we assume that $Z$ is chosen to be the minimal closed
subset satisfying equation~(\ref{conjeqn}).

Let $V$ be an irreducible component of $Z$.  We will show that $V$
does not surject onto $X$ via both projections from $Y$ unless $V=D$.
Thus, assume that $V$ does surject onto $X$ via both projections.
Then $\dim V\geq\dim X$.  Let $E=V\cap D$.  If $E=V$, then $V=D$, so
we may assume that $E$ is a cycle with positive codimension on $V$.
By, for example, Theorems 3.26 and 3.27 from \cite{Ko}, there is a
smooth variety $\tilde{V}$ and a dominant birational map
$\pi\colon\tilde{V}\to V$ such that $\pi$ is an isomorphism away from
$N=E\cup\mbox{Sing}(V)$, and such that the induced reduced cycle $M$
of $\pi^*N$ is contained in a divisor with normal crossings.

If $\tilde{V}$ has negative Kodaira dimension, then by assumption it
is uniruled.  Since it is birational to a subvariety of $X\times X$
that surjects onto $X$ via both projections, this implies that $X$ is
also uniruled.  Since $X$ has nonnegative Kodaira dimension, this is
impossible.

Thus, we may assume that $\tilde{V}$ has nonnegative Kodaira
dimension.  This means that some multiple of the canonical divisor of
$\tilde{V}$ must admit a global section, and so the height $h_K$ is
bounded below by a constant.  In addition, if $M$ is the sum of the
irreducible components of $\pi^*N$ (i.e., $M$ is the reduced divisor
induced by $\pi^*N$), then:
\[m_S(\pi^*N,P)+h_K(P)\leq \alpha(m_S(M,P)+h_K(P))+O(1)\]
for some positive integer $\alpha$.  Thus, if Vojta's Main Conjecture
(Conjecture~\ref{vojta}) is true for $M$, it must also be true for
$\pi^*N$.  (Note that $M$ is, as noted above, contained in a normal
crossings divisor.)

Thus, by Vojta's Main Conjecture, we find a proper Zariski closed
subset $\tilde{Z}$ of $\tilde{V}$ such that for all $P\in
\tilde{V}(k)-\tilde{Z}(k)$, we have:
\[m_{E,v}(P)\leq\epsilon h(P) + O(1)\]
where $h(P)$ in this case denotes $h(\pi(P))$, since this yields a big
height on $\tilde{V}$.  Noting that $m_{D,v}(\pi(P))= m_{E,v}(P)+O(1)$
for $P\in\tilde{V}(k)=E(k)$ gives:
\[m_{D,v}(P,Q)\leq\epsilon h(P,Q) + O(1)\]
for all $k$-rational points $(P,Q)$ in some nonempty open subset of
$V$.  But this is precisely equation~(\ref{conjeqn}), which
contradicts the minimality of $V$.  We conclude that $V$ does not
surject onto $X$ via both projections unless $V=D$.

Let $W$ be the union of all proper closed subsets $H$ of $X$ such that
$H$ is the projection of some irreducible component of $Z$.  (Recall
that $Z$ is the exceptional subset of $X$ derived from Vojta's Main
Conjecture.)  Let $U'=X-W$.  Then the intersection of $U'\times U'$
with the exceptional set $Z\subset X\times X$ is a subset of the
diagonal, so for each $P,Q\in U'(k)$ with $P\neq Q$, we have
$m_{D,v}(P,Q) < \epsilon h(P,Q) + c$.  Taking the reciprocal
exponential of both sides of that equation yields:
\[\dist_v(P,Q) > CB^{-\epsilon}\]
as desired.  \qed

\vspace{.1in}

Note that the Repulsion Principle is true unconditionally for
subvarieties of abelian varieties.  This is because Vojta's Main
Conjecture for subvarieties of abelian varieties was proven by
Faltings (\cite{Fa2}).  (Abelian varieties do not contain any
subvarieties of negative Kodaira dimension.)

\section{Vojta's Conjecture implies the Batyrev-Manin Conjecture for 
$K3$ surfaces}\label{bm}

Before we can prove the implication in the title of this section, we
must describe the relevant conjecture of Batyrev and Manin.  For
context and motivation of this conjecture, we refer the reader to
\cite{BM}.

Let $X$ be an algebraic variety defined over a number field $k$, and
let $L$ be an ample divisor on $X$.  Choose a (multiplicative) height
function $H_L$ on $V$ corresponding to $L$, and let $W$ be any subset
of $X$.  We define the counting function for $W$ by:
\[N_{W,L}(B) = \#\{P\in S(k)\mid H_L(P)<B\}\]
Batyrev and Manin \cite{BM} have made a series of conjectures about the
behaviour of $N_{W,L}$.  In the case that $X$ is a $K3$ surface, their
conjecture is as follows:

\begin{conj}[Batyrev-Manin Conjecture for $K3$ Surfaces]\label{batyrevmanin}
Let $\epsilon>0$ be any real number and $L$ any ample divisor on $X$.
Then there is a non-empty Zariski open subset $U(\epsilon)\subset X$
such that
\[N_{U(\epsilon),L}(B) = O(B^\epsilon)\]
\end{conj}

In fact, this principle should apply much more broadly than just to $K3$
surfaces.  We have the following general result:

\begin{thm}\label{main}
Let $X$ be a smooth, projective variety of nonnegative Kodaira
dimension, defined over a number field $k$.  Let $V$ be a smooth
variety birational to a subvariety of $X\times X$, and make the
following two assumptions:
\begin{enumerate}
\item If $V$ has nonnegative Kodaira dimension, assume that Vojta's Main
Conjecture for general cycles (Conjecture~\ref{vojta}) is true for $V$.
\item If $V$ has negative Kodaira dimension, assume that $V$ is uniruled.
\end{enumerate}
Then for every $\epsilon>0$ and ample divisor $L$ on $X$, there is a
nonempty Zariski open subset $U(\epsilon)\subset X$ such that
\[N_{U(\epsilon),L}(B) = O(B^\epsilon)\]
\end{thm}

\noindent
{\it Proof:} \/ Let $v$ be an archimedean place of $k$, and let $n$ be
the real dimension of the variety $X_v$ over the associated completion
$k_v$ of $k$.  (That is, $n=\dim X$ if $v$ is a real place, and
$n=2\dim X$ if $v$ is a complex place.)  By the Repulsion Principle,
there is a nonempty Zariski open subset $U$ of $X$ and a positive real
constant $C$ such that for all $P$ and $Q$ in $U(k)$, we have:
\[\dist_v(P,Q) > CB^{-\epsilon/n}\]
for all sufficiently large real numbers $B$.

Let $A$ be the set $A=\{P\in U(k)\mid H(P)\leq B\}$.
For any pair of points $P$ and $Q$ in $A$, we have:
\[\dist_v(P,Q)> CB^{-\epsilon/n}\]
Since $A_i$ lies in a real manifold of finite dimension $n$ and finite
volume, it follows that:
\[\# A=O(B^{\epsilon})\]
as desired.  \qed

\vspace{.1in}

Mori has proven (\cite{Mo}) that any smooth projective threefold with
negative Kodaira dimension is uniruled.  Since curves and surfaces
with negative Kodaira dimension are also well known to be uniruled, we
deduce the following corollary of Theorem~\ref{main}:

\begin{cor}\label{vojtabm}
Let $X$ be a smooth, projective surface of nonnegative Kodaira
dimension, defined over a number field $k$.  Assume that Vojta's Main
Conjecture for general cycles (Conjecture~\ref{vojta}) is true for any
variety birational to a subvariety of $X\times X$.  Then for every
$\epsilon>0$ and ample divisor $L$ on $X$, there is a nonempty Zariski
open subset $U(\epsilon)\subset X$ such that
\[N_{U(\epsilon),L}(B) = O(B^\epsilon)\]
In particular, Vojta's Main Conjecture implies the Batyrev-Manin
Conjecture for $K3$ surfaces.
\end{cor}

\section{Manin's Rational Curve Conjecture}\label{rcc}

One can almost generalize the arguments here to prove the Rational
Curve Conjecture of Manin (\cite{Ma}), if one further assumes a
conjecture implied by the Minimal Model Program.  Here is Manin's 
conjecture:

\begin{conj}[Rational Curve Conjecture]\label{rcconj}
Let $U$ be a nonempty Zariski open subset of a smooth algebraic
variety $X$ defined over a number field $k$, and let $L$ be an ample
divisor on $X$.  If there are positive real constants $\delta$ and $c$
such that $N_{U,L}(B)\geq cB^\delta$ for infinitely many arbitrarily
large positive real numbers $B$, then $U$ contains a nonempty open
subset of a rational curve $C$, defined over $k$ and containing a
dense set of $k$-rational points.
\end{conj}

We will not quite be able to prove this.  Instead, we will prove
(modulo two conjectures) the slightly weaker assertion that $U$ must
contain a nonempty open subset of a rational curve defined over some
finite extension of $k$.

\begin{thm}\label{vojtarcc}
Let $U$ be a nonempty Zariski open subset of a smooth algebraic
variety $X$ defined over a number field $k$, and let $L$ be an ample
divisor on $X$.  For any smooth variety $V$ birational to a subvariety of
$X\times X$, make the following two assumptions:
\begin{enumerate}
\item If $V$ has nonnegative Kodaira dimension, assume that Vojta's Main
Conjecture for general cycles (Conjecture~\ref{vojta}) is true for $V$.
\item If $V$ has negative Kodaira dimension, assume that $V$ is uniruled.
\end{enumerate}
If there are positive real constants $\delta$ and $c$ such that
$N_{U,L}(B)\geq cB^\delta$ for infinitely many arbitrarily large
positive real numbers $B$, then $U$ contains a nonempty open subset of
a rational curve $C$ (not necessarily defined over $k$).
\end{thm}

\noindent
{\it Proof:} \/ First, note that if the Kodaira dimension of $X$ is
negative, then the result is trivially true.  Thus, in what follows,
we assume that $X$ has nonnegative Kodaira dimension.  In particular,
we may assume that some multiple of the canonical class admits a
global section, and therefore that any height $h_K$ associated to the
canonical class is bounded from below by a constant.

We will induce on the dimension of $X$.  Theorem~\ref{vojtarcc} is
clearly true if $\dim X=1$, so assume that it is true for all
varieties of dimension less than $X$.  By Theorem~\ref{main}, there is
a nonempty Zariski open subset $U\subset X$ such that
\[N_{U,L}(B)=O(B^\epsilon)\]
By comparing counting functions, this implies for small enough
$\epsilon$ that there is some proper closed subset $Y$ of $U$ which is
not contained in $U$.  Moreover, we must have $\#\{P\in Y(k)\mid
H(P)\leq B\} > CB^\delta$ for infinitely many arbitrarily large
positive real numbers $B$ and some real positive constant $C$.  By
the induction hypothesis, $Y\subset U$ must contain a nonempty open
subset of a rational curve, as desired.  \qed

\vspace{.1in}

As noted previously, the assumption that every variety of negative
Kodaira dimension is uniruled is known to be true in dimension at most
three, so the Rational Curve Conjecture for surfaces follows simply
from Vojta's Main Conjecture.  In higher dimensions, the standard
conjectures of the Minimal Model Program imply that every variety of
negative Kodaira dimension is uniruled.  The Hard Dichotomy Theorem
implies that the result of the Minimal Model Program applied to any
variety of negative Kodaira dimension is a Mori fibre space.  It is a
well known result of Miyaoka and Mori (\cite{MM}) that any Mori fibre
space is uniruled.  For details, see for example \cite{Mat}, section
3.2.

\end{document}